\documentclass[10pt,reqno]{amsart}
\usepackage[OT1]{fontenc}
\usepackage{amsmath,amssymb,amsthm}
\usepackage[a4paper,hcentering,textwidth=138mm,top=25mm,bottom=24mm,
  headheight=10pt,headsep=12pt,footskip=20pt]{geometry}
\usepackage{etoolbox}
\usepackage[hidelinks]{hyperref}
\hypersetup{pdftitle={The inverse Erdos-Fuchs problem for k-fold sums}}
\newtheorem{priorresult}{Theorem}

\newtheorem{theorem}{Theorem}[section]
\newtheorem{lemma}[theorem]{Lemma}
\newtheorem{proposition}[theorem]{Proposition}
\numberwithin{equation}{section}
\newcommand{\N}{\mathbb{N}}
\newcommand{\Nzero}{\mathbb{N}_0}
\newcommand{\E}{\mathbb{E}}
\newcommand{\Prob}{\mathbb{P}}
\DeclareMathOperator{\vol}{vol}

\title[Inverse Erd\H{o}s--Fuchs bounds for $k$-fold sums]
{The inverse Erd\H{o}s--Fuchs problem for $k$-fold sums}
\author[Aimin Guo]{Aimin Guo}
\author[Xiaohui Yan]{Xiaohui Yan}
\thanks{Corresponding author: Xiaohui Yan
  (\texttt{yanxiaohui\_1992@163.com}).}
\address[Aimin Guo and Xiaohui Yan]{School of Mathematics and Statistics,
  Anhui Normal University\newline
  \indent Wuhu 241002, Anhui, P.~R.~China}
\email[Aimin Guo]{18339681864@163.com}
\email[Xiaohui Yan]{yanxiaohui\_1992@163.com}
\date{November 11, 2026}
\makeatletter
\renewcommand{\section}{\@startsection{section}{1}{\z@}%
  {12pt plus 2pt minus 1pt}{7pt}{\normalfont\bfseries\centering}}
\patchcmd{\@maketitle}{\topskip42\p@}{\topskip30\p@}{}{}
\patchcmd{\@maketitle}{\dimen@34\p@}{\dimen@19\p@}{}{}
\@ifundefined{subjclassname@2020}{%
  \@namedef{subjclassname@2020}{\textup{2020} Mathematics Subject Classification}%
}{}
\makeatother

\begin{document}
\begin{abstract}
For every fixed integer \(k\ge3\), we construct an infinite set
\(A\subseteq\N\) such that
\[
 \sum_{h\le n}r_{kA}(h)
 =\Gamma(1+1/k)^k n+O_k\bigl(n^{(1-1/k)^2}\bigr),
\]
where \(r_{kA}(h)\) counts ordered representations of \(h\) as a
sum of \(k\) elements of \(A\), with repetitions allowed.
This improves Dai and Pan's bound for a linear main term.
We also prove that, almost surely in this paired random model,
the absolute discrepancy is at least a positive constant multiple of
\(n^{(1-1/k)^2}\) for infinitely many \(n\).
\end{abstract}
\maketitle

\section{Introduction}

Let \(\N=\{1,2,\ldots\}\) and \(\Nzero=\{0,1,2,\ldots\}\).
For an integer \(k\ge2\), a set \(A\subseteq\Nzero\), and
\(h\in\Nzero\), define
\[
 r_{kA}(h)=\#\{(b_1,\ldots,b_k)\in A^k:b_1+\cdots+b_k=h\}.
\]
The tuples are ordered and repetitions are allowed.  Thus
\[
 \sum_{h\le n}r_{kA}(h)
 =\#\{(b_1,\ldots,b_k)\in A^k:b_1+\cdots+b_k\le n\},
\]
where \(h\) runs over the nonnegative integers.  We first recall
the theorem of Erd\H{o}s and Fuchs.

\begin{priorresult}[\cite{EF}]
Let \(A\subseteq\Nzero\) be infinite and let \(c>0\).  Then
\[
 \sum_{h\le n}r_{2A}(h)=cn+o\bigl(n^{1/4}(\log n)^{-1/2}\bigr)
\]
cannot hold as \(n\to\infty\).
\end{priorresult}

Jurkat, in unpublished work, and
Montgomery and Vaughan \cite{MV} independently ruled out even
\(o(n^{1/4})\).  Tang \cite{Tang} studied lower bounds for sums
with more terms.  For related work on the optimal exponent in the
Erd\H{o}s--Fuchs theorem for \(k\)-fold sums, see Yan \cite{Yan}.

In the converse direction, Ruzsa constructed the following set.

\begin{priorresult}[\cite{Ruzsa}]
There is an infinite set \(A\subseteq\Nzero\) such that
\[
 \sum_{h\le n}r_{2A}(h)=\frac{\pi}{4}n+O(n^{1/4}\log n).
\]
\end{priorresult}

Yang and Zhao improved the logarithmic factor.

\begin{priorresult}[{\cite[Theorem~1.1]{YZ}}]\label{thm:yang-zhao}
There is an infinite set \(A\subseteq\Nzero\) such that
\[
 \sum_{h\le n}r_{2A}(h)=\frac{\pi}{4}n+O(n^{1/4}\sqrt{\log n}).
\]
\end{priorresult}

Their proof pairs adjacent random points and uses the bounded
differences inequality \cite{McD}.  We use the same pairing for
\(k\ge3\).

Dai and Pan \cite{DP} studied \(k\)-fold sums with main term \(Cn^\beta\).
Taking \(\beta=1\) in their Theorem~1.1 yields the following bound.

\begin{priorresult}[{\cite[Theorem~1.1]{DP}}]
For each fixed integer \(k\ge3\), there is an infinite set
\(A\subseteq\Nzero\) and a constant \(C>0\) such that
\begin{equation}\label{eq:old}
 \sum_{h\le n}r_{kA}(h)=Cn+O_k\bigl(n^{1-(k+1)/k^2}\sqrt{\log n}\bigr).
\end{equation}
\end{priorresult}

The constant \(C\) is left unevaluated in their theorem statement.
For \(\beta=1\), their construction identifies it as the volume of
\(\{y\in[0,\infty)^k:\sum_\nu y_\nu^k\le1\}\).
Evaluating this volume by \eqref{eq:volume} gives
\(C=\Gamma(1+1/k)^k\).  We retain this main term and improve the
error term as follows.

\begin{theorem}\label{thm:main}
For every fixed integer \(k\ge3\), there is an infinite set
\(A=A_k\subseteq\N\) such that, as the integer \(n\) tends to infinity,
\begin{equation}\label{eq:main}
 \sum_{h\le n}r_{kA}(h)=\Gamma\!\left(1+\frac1k\right)^kn
            +O_k\bigl(n^{1-2/k+1/k^2}\bigr).
\end{equation}
Moreover, for real \(X\ge1\),
\[
 A(X):=\#(A\cap[1,X])=X^{1/k}+O(1).
\]
\end{theorem}

The exponent in \eqref{eq:main} is smaller than that in \eqref{eq:old}
by \((k-2)/k^2\), and the logarithmic factor is removed.  In particular,
for \(k=3\), we obtain \(O(n^{4/9})\) in place of
\(O(n^{5/9}\sqrt{\log n})\).

Throughout the proof, \(k\ge3\) is fixed and \(u_+=\max(u,0)\).
We write \(\E\) for expectation and \(\Prob\) for probability.
In \(O_k\), \(\ll_k\), and \(\asymp_k\), the subscript allows
constants depending on \(k\); \(\asymp_k\) means \(\ll_k\) in both
directions.

\section{The random construction}\label{sec:construction}

Take independent random variables \(U_0,U_1,\ldots\), each uniform on
\((0,1)\).  For each pair index \(m\in\Nzero\), put
\begin{equation}\label{eq:points}
 \theta_{2m}=2m+U_m,\qquad \theta_{2m+1}=2m+2-U_m .
\end{equation}
Each sample point \(\theta_s\) is uniform on \((s,s+1)\).
We use \(s\) for a sample index and \(m\) for a pair index.
The interval
\(I_m=[2m,2m+2)\) contains the two points determined by \(U_m\),
one in each unit interval; different pairs are independent.

To obtain a set of positive integers, define
\[
 \begin{aligned}
 a_s&=1+\lfloor\theta_s^k\rfloor,\qquad A=\{a_s:s\ge0\},\\
 s^k+1&\le a_s\le(s+1)^k .
 \end{aligned}
\]
The bounds follow from \(s<\theta_s<s+1\).  These integer ranges
are disjoint, so each element of \(A\) has a unique sample index \(s\).

Fix an integer \(n\ge1\) and put \(R=n^{1/k}\).  For ordered sample
indices \(s_1,\ldots,s_k\ge0\), with repetitions allowed, define
the count before rounding by
\[
 \begin{aligned}
 \delta_{s_1,\ldots,s_k}
 &=\mathbf{1}_{\{\theta_{s_1}^k+\cdots+\theta_{s_k}^k\le n\}},\\
 \sigma_n&=\sum_{s_1,\ldots,s_k\ge0}\delta_{s_1,\ldots,s_k}.
 \end{aligned}
\]
Here \(\mathbf{1}_E\) denotes the indicator of \(E\).  Only indices
\(s_\nu<R\) can contribute, so the sums are finite.
The bounds \(\theta_s^k\le a_s\le\theta_s^k+1\) give
\begin{equation}\label{eq:rounding}
 \begin{gathered}
 \sum_{\nu=1}^k\theta_{s_\nu}^k
 \le\sum_{\nu=1}^k a_{s_\nu}
 \le\sum_{\nu=1}^k\theta_{s_\nu}^k+k,\\
 \sigma_{n-k}\le\sum_{h\le n}r_{kA}(h)\le\sigma_n
 \qquad(n\ge k+1).
 \end{gathered}
\end{equation}
It therefore suffices to compare \(\sigma_n\) with the volume of
\(\Omega_n=\{y\in[0,\infty)^k:y_1^k+\cdots+y_k^k\le n\}\).
We will also need slice volumes.  With \(\vol\) denoting volume in
the ambient dimension, for an integer \(d\ge1\) and real \(t\ge0\)
we have
\begin{equation}\label{eq:volume}
 \vol\{y\in[0,\infty)^d:y_1^k+\cdots+y_d^k\le t\}
 =C_{d,k}t^{d/k},\qquad
 C_{d,k}=\frac{\Gamma(1+1/k)^d}{\Gamma(1+d/k)}.
\end{equation}
For \(t>0\), substitute \(u_\nu=y_\nu^k/t\) in each coordinate.
Then \(dy_\nu=(t^{1/k}/k)u_\nu^{1/k-1}\,du_\nu\), and the region
becomes the simplex \(\sum_\nu u_\nu\le1\), \(u_\nu\ge0\).
The volume is
\[
 \frac{t^{d/k}}{k^d}
 \int_{\substack{u_\nu\ge0\\ \sum_\nu u_\nu\le1}}
       \prod_{\nu=1}^d u_\nu^{1/k-1}\,du_1\cdots du_d
 =
 \frac{t^{d/k}}{k^d}\frac{\Gamma(1/k)^d}{\Gamma(1+d/k)} .
\]
The equality is the Dirichlet integral on the simplex.  Since
\(\Gamma(1+1/k)=\Gamma(1/k)/k\), this gives \eqref{eq:volume}.
For \(t=0\), both sides are zero.
For \(d=k\) and \(t=n\), this is \(C_kn\), where
\(C_k=C_{k,k}=\Gamma(1+1/k)^k\).

We now partition the nonnegative orthant into boxes
\(B=I_{m_1}\times\cdots\times I_{m_k}\), where \(m_\nu\in\Nzero\)
is the pair index in coordinate \(\nu\).  The count and volume in
each box are
\[
 \delta_B=\sum_{\substack{s_\nu\in\{2m_\nu,2m_\nu+1\}\\1\le\nu\le k}}
                 \delta_{s_1,\ldots,s_k},
 \qquad d_B=\vol(B\cap\Omega_n).
\]
Two coordinates with the same pair index have the same two choices
of sample point.  Summing over the boxes gives
\[
 \sigma_n-C_kn
 =\sum_B\delta_B-\sum_Bd_B
 =\sum_B(\delta_B-d_B).
\]
For boxes entirely inside or outside \(\Omega_n\), respectively,
\[
 \delta_B=d_B=
 \begin{cases}
  2^k,&\displaystyle\sum_{\nu=1}^k(2m_\nu+2)^k\le n,\\
  0,&\displaystyle\sum_{\nu=1}^k(2m_\nu)^k>n.
 \end{cases}
\]

By coordinatewise monotonicity, only boxes satisfying
\[
 \sum_{\nu=1}^k(2m_\nu)^k\le n
 <\sum_{\nu=1}^k(2m_\nu+2)^k
\]
can contribute; we call them boundary boxes.
For distinct pair indices, averaging the count recovers the volume.
Indeed, for any bounded measurable \(g\) on \(I_m\),
\begin{equation}\label{eq:pair-uniform}
 \begin{aligned}
 \E\bigl[g(\theta_{2m})+g(\theta_{2m+1})\bigr]
 &=\int_0^1\bigl(g(2m+u)+g(2m+2-u)\bigr)\,du\\
 &=\int_{2m}^{2m+1}g(y)\,dy+\int_{2m+1}^{2m+2}g(y)\,dy\\
 &=\int_{I_m}g(y)\,dy .
 \end{aligned}
\end{equation}
Apply this identity to
\(g(y)=\mathbf{1}_{\{y^k+y_2^k+\cdots+y_k^k\le n\}}\).
When \(m_1,\ldots,m_k\) are distinct, independence allows us to
average successively over the coordinates, giving
\[
 \E\delta_B
 =\int_{I_{m_k}}\cdots\int_{I_{m_1}}
 \mathbf{1}_{\{y_1^k+\cdots+y_k^k\le n\}}\,dy_1\cdots dy_k
 =d_B .
\]
Let \(\sigma_n^\ast\) count the tuples from boxes with distinct
pair indices.

It remains to account for boxes with repeated pair indices, whose
coordinates can use the same draw.  Since \(0\le\delta_B,d_B\le2^k\),
we have \(|\delta_B-d_B|\le2^k\); the difference vanishes off the
boundary.  It therefore suffices to count these boundary boxes.

\begin{lemma}\label{lem:boundary}
There are \(O_k(R^{k-2})\) ordered boundary boxes with a repeated
pair index.
\end{lemma}

\begin{proof}
Choose the two equal coordinate positions in \(\binom{k}{2}\) ways.
Writing them first, the boxes have the form
\[
 I_{z_1}\times I_{z_1}\times I_{z_2}\times\cdots\times I_{z_{k-1}},
\]
where \(z=(z_1,\ldots,z_{k-1})\in\Nzero^{k-1}\); further
repetitions are allowed.  The lower corner sum is
\[
 \Phi(z)=2(2z_1)^k+\sum_{\nu=2}^{k-1}(2z_\nu)^k .
\]
With \(\mathbf{1}=(1,\ldots,1)\), the boundary condition is
\(\Phi(z)\le n<\Phi(z+\mathbf{1})\).

Subtract the smallest entry of \(z\) from each entry.  This gives
a unique start \(w\) and a nonnegative integer \(\lambda\) such that
\[
 w=z-\lambda\mathbf{1},\qquad
 \lambda=\min_\nu z_\nu,\qquad \min_\nu w_\nu=0.
\]
The list from \(w\) consists of \(w+\lambda\mathbf{1}\), \(\lambda\ge0\).
For \(k=3\), the list containing \((3,5)\) is
\((0,2),(1,3),(2,4),(3,5),\ldots\).
The supremum of the power sum on one box is the minimum on the next.
Consequently, the intervals
\[
 \bigl[\Phi(w+\lambda\mathbf{1}),\,
       \Phi(w+(\lambda+1)\mathbf{1})\bigr),\qquad\lambda\ge0,
\]
are disjoint.  Each list therefore contains at most one boundary
box for a fixed \(n\).

For a boundary box, \(0\le w_\nu\le z_\nu\le\lfloor R/2\rfloor\).
Choose a zero position in \(w\) in at most \(k-1\) ways and each
remaining entry in at most \(\lfloor R/2\rfloor+1\) ways.
If several entries of \(w\) vanish, a start may be counted more
than once; this is harmless for an upper bound.
Every repeated-index box occurs for some choice of the two equal
positions, so the number of boxes is at most
\[
 \binom{k}{2}(k-1)\bigl(\lfloor R/2\rfloor+1\bigr)^{k-2}
 \ll_k R^{k-2}.\qedhere
\]
\end{proof}

Since \(\E\delta_B=d_B\) for distinct indices and
\(|\delta_B-d_B|\le2^k\) otherwise, Lemma~\ref{lem:boundary} gives
\begin{equation}\label{eq:reduce}
 \begin{aligned}[b]
  \sigma_n-C_kn
  &=\sum_{\substack{B\text{ with distinct}\\\text{pair indices}}}(\delta_B-d_B)
    +\sum_{\substack{B\text{ boundary with}\\\text{a repeated pair index}}}(\delta_B-d_B)\\
  &=\sigma_n^\ast-\E\sigma_n^\ast+O_k(R^{k-2}).
 \end{aligned}
\end{equation}

We record three probability tools for \eqref{eq:reduce}.
The Hoeffding decomposition separates the draws, a martingale
moment bound controls their interactions, and the first
Borel--Cantelli lemma converts summable tail probabilities into
an almost sure estimate.
For a random variable \(Y\), write \(\|Y\|_p=(\E|Y|^p)^{1/p}\)
for \(1\le p<\infty\).

\begin{lemma}[Hoeffding decomposition \cite{Hoeffding}]\label{lem:hoeffding}
Let \(S\) be a finite index set, let \((U_m)_{m\in S}\) be
independent random variables, and let \(W\) be a bounded measurable
real-valued function of these variables.  For \(Q\subseteq S\),
write \(U_Q=(U_m)_{m\in Q}\).  For \(T\subseteq S\), put
\[
 W_T=\sum_{K\subseteq T}(-1)^{|T|-|K|}\E[W\mid U_K],
\]
where \(|T|\) is the number of elements of \(T\), and conditioning
on \(U_\varnothing\) means taking expectation.  Then
\[
 W_\varnothing=\E W,\qquad
 W-\E W=\sum_{\varnothing\ne T\subseteq S}W_T.
\]
Each nonempty component has zero conditional mean in any one
of its variables, and conditioning on a fixed set \(Q\) gives
\begin{equation}\label{eq:hoeffding-conditioning}
 \begin{aligned}
 \E[W_T\mid U_{T\setminus\{m\}}]&=0
       &&(\varnothing\ne T\subseteq S,\ m\in T),\\
 \E[W\mid U_Q]&=\sum_{T\subseteq Q}W_T
       &&(Q\subseteq S),
 \end{aligned}
\end{equation}
almost surely.
\end{lemma}

\begin{proof}
Fix \(m\in T\).  Pair each \(K\subseteq T\setminus\{m\}\)
with \(K\cup\{m\}\).  By independence, averaging over \(U_m\)
in \(\E[W\mid U_{K\cup\{m\}}]\) gives \(\E[W\mid U_K]\).
The paired terms have opposite signs, proving the conditional
mean-zero assertion.

For fixed \(K\subseteq S\), the coefficient of \(\E[W\mid U_K]\)
in \(\sum_{T\subseteq S}W_T\) is
\[
 \sum_{T:\,K\subseteq T\subseteq S}(-1)^{|T|-|K|}
 =\begin{cases}1,&K=S,\\0,&K\ne S.\end{cases}
\]
Thus \(\sum_TW_T=W\), giving the decomposition.  If
\(T\not\subseteq Q\), independence and the conditional mean-zero
property give \(\E[W_T\mid U_Q]=0\); otherwise this expectation is
\(W_T\).  Conditioning the decomposition proves
\eqref{eq:hoeffding-conditioning}.
\end{proof}

\begin{lemma}[Moment bound for martingale differences]\label{lem:martingale}
Let \(\mathcal F_0\subseteq\cdots\subseteq\mathcal F_M\) be
an increasing sequence of sigma algebras.  Suppose that \(Z_\nu\)
is bounded and real-valued, is measurable with respect to
\(\mathcal F_\nu\), and satisfies
\(\E[Z_\nu\mid\mathcal F_{\nu-1}]=0\) for \(1\le\nu\le M\).
Then, for \(2\le p<\infty\),
\begin{equation}\label{eq:martingale}
 \left\|\sum_{\nu=1}^M Z_\nu\right\|_p^2
 \le(p-1)\sum_{\nu=1}^M\|Z_\nu\|_p^2.
\end{equation}
\end{lemma}

\begin{proof}
For a fixed step \(\nu\), put \(V=\sum_{a=1}^{\nu-1}Z_a\) and
\(Y=Z_\nu\).
Conditional centering gives
\(\E[VY]=\E[V\E[Y\mid\mathcal F_{\nu-1}]]=0\), proving the case
\(p=2\) through
\(\|V+Y\|_2^2=\|V\|_2^2+\|Y\|_2^2\).

For \(p>2\) and \(\varepsilon>0\), define, within this proof,
\(f(\tau)=\varepsilon+\E|V+\tau Y|^p\) and
\(g(\tau)=f(\tau)^{2/p}\) for \(0\le\tau\le1\).
The boundedness of \(V,Y\) allows differentiation under the
expectation.  As \(|V|^{p-2}V\) is
\(\mathcal F_{\nu-1}\)-measurable, we have
\[
 f'(0)=p\E[|V|^{p-2}VY]
      =p\E\bigl[|V|^{p-2}V\,\E[Y\mid\mathcal F_{\nu-1}]\bigr]=0.
\]
Hence \(g'(0)=\frac2p f(0)^{2/p-1}f'(0)=0\).
H\"older's inequality with exponents \(p/(p-2)\) and \(p/2\),
followed by the chain rule, gives
\[
 \begin{aligned}
 f''(\tau)
 &=p(p-1)\E[|V+\tau Y|^{p-2}Y^2]\\
 &\le p(p-1)\|V+\tau Y\|_p^{p-2}\|Y\|_p^2
 \le p(p-1)f(\tau)^{1-2/p}\|Y\|_p^2,\\
 g''(\tau)
 &=\frac2p\left(\frac2p-1\right)f(\tau)^{2/p-2}
       \bigl(f'(\tau)\bigr)^2
   +\frac2p f(\tau)^{2/p-1}f''(\tau)\\
 &\le \frac2p f(\tau)^{2/p-1}f''(\tau)
  \le 2(p-1)\|Y\|_p^2 .
 \end{aligned}
\]
In the estimate for \(g''\), the squared-derivative term is
nonpositive because \(2/p-1<0\).  Consequently,
\[
 g(1)-g(0)
 =g'(0)+\int_0^1(1-\tau)g''(\tau)\,d\tau
 \le(p-1)\|Y\|_p^2.
\]
Letting \(\varepsilon\downarrow0\) gives
\(\|V+Y\|_p^2-\|V\|_p^2\le(p-1)\|Y\|_p^2\).
Summing these one-step inequalities for \(\nu=1,\ldots,M\) proves
\eqref{eq:martingale}.
\end{proof}

\begin{lemma}[First Borel--Cantelli lemma {\cite[Theorem~2.3.1]{Durrett}}]\label{lem:borel-cantelli}
Let \((\mathcal E_n)_{n\ge2}\) be measurable events on a common
probability space.  If
\[
 \sum_{n=2}^{\infty}\Prob(\mathcal E_n)<\infty,
\]
then, with probability one, only finitely many of the events
\(\mathcal E_n\) occur.  Equivalently,
\[
 \Prob\left(\bigcap_{N=2}^{\infty}
                \bigcup_{n\ge N}\mathcal E_n\right)=0.
\]
Independence of the events is not required.
\end{lemma}

\begin{proof}
For each integer \(N\ge2\), the union bound gives
\[
 \Prob\left(\bigcup_{n\ge N}\mathcal E_n\right)
 \le\sum_{n\ge N}\Prob(\mathcal E_n)
 \longrightarrow0\qquad(N\to\infty).
\]
If infinitely many \(\mathcal E_n\) occur, then at least one occurs
after every starting index \(N\).  This event is contained in each
tail union above, so its probability is zero.  Outside this
probability-zero event, there is a finite starting index after
which no \(\mathcal E_n\) occurs.
\end{proof}

\section{Proof of Theorem~\ref{thm:main}}\label{sec:proof}

By \eqref{eq:reduce}, it remains to estimate
\(\sigma_n^\ast-\E\sigma_n^\ast\).  Distinct-index boxes can share
draws, so we separate individual draws from their interactions.
Only pair indices in \(S=\{m\in\Nzero:2m\le R\}\) can contribute.
As in Lemma~\ref{lem:hoeffding}, \(U_T\) denotes the draws indexed
by \(T\); conditioning on \(U_T\) fixes these draws and averages
over the others.  The effect of a single draw is
\[
 G_{\{m\}}=\E[\sigma_n^\ast\mid U_m]-\E\sigma_n^\ast .
\]
For two distinct pairs \(m_1,m_2\), subtract the change due to each
pair alone from the change due to fixing both:
\[
 \begin{aligned}
 G_{\{m_1,m_2\}}
 ={}&\bigl(\E[\sigma_n^\ast\mid U_{m_1},U_{m_2}]
              -\E\sigma_n^\ast\bigr)\\
 &-\bigl(\E[\sigma_n^\ast\mid U_{m_1}]-\E\sigma_n^\ast\bigr)
  -\bigl(\E[\sigma_n^\ast\mid U_{m_2}]-\E\sigma_n^\ast\bigr)\\
 ={}&\E[\sigma_n^\ast\mid U_{m_1},U_{m_2}]
  -\E[\sigma_n^\ast\mid U_{m_1}]\\
 &-\E[\sigma_n^\ast\mid U_{m_2}]+\E\sigma_n^\ast .
 \end{aligned}
\]
More generally, for \(T\subseteq S\), define
\[
 G_T=\sum_{K\subseteq T}(-1)^{|T|-|K|}\E[\sigma_n^\ast\mid U_K].
\]
Here \(K\) is the summation subset; for \(K=\varnothing\), the
conditional expectation means \(\E\sigma_n^\ast\), so
\(G_\varnothing=\E\sigma_n^\ast\).
We suppress the dependence of \(G_T\) on \(n\) in the notation.
Lemma~\ref{lem:hoeffding}, applied to \(W=\sigma_n^\ast\), gives
\[
 \sum_{\varnothing\ne T\subseteq S}G_T=\sigma_n^\ast-\E\sigma_n^\ast,
 \qquad \E[G_T\mid U_{T\setminus\{m\}}]=0\quad(m\in T).
\]

If a box omits an index in \(T\), its contribution to \(G_T\)
cancels: fixing or averaging over that index gives the same count
in the box.  In particular, \(G_T=0\) for \(|T|>k\), since each
box uses at most \(k\) pairs.
We now separate the single-pair terms from all the interactions:
\begin{equation}\label{eq:decomp}
 \begin{aligned}
 \sigma_n^\ast-\E\sigma_n^\ast
 &=\sum_{m\in S}G_{\{m\}}+\mathcal H_n,\\
 \mathcal H_n
 &=\sum_{\substack{T\subseteq S\\2\le |T|\le k}}G_T.
 \end{aligned}
\end{equation}
In estimates for these terms, \(m\) ranges over \(S\) and \(T\)
over subsets of \(S\).

Combining \eqref{eq:reduce} and \eqref{eq:decomp}, we obtain
\begin{equation}\label{eq:total-discrepancy}
 \sigma_n-C_kn=\sum_m G_{\{m\}}+\mathcal H_n+O_k(R^{k-2}).
\end{equation}
Put \(\gamma=k-2+1/k\), so that \(R^\gamma=n^{(1-1/k)^2}\).
We first bound the sum of the single-pair terms by \(O_k(R^\gamma)\)
for every choice of the points.  We then show that
\(\mathcal H_n=o(R^\gamma)\) almost surely.

For the individual draws, we use the volumes of slices with one
coordinate fixed.  By \eqref{eq:volume}, the slice at \(x\ge0\)
and the sum of the two slices from pair \(m\) are
\[
 f_n(x)=C_{k-1,k}(n-x^k)_+^{1-1/k},\qquad
 q_{m,n}(u)=f_n(2m+u)+f_n(2m+2-u).
\]
Here \(f_n(x)=0\) for \(x\ge R\), and \(0<u<1\).
In \(q_{m,n}\), the first subscript is the pair index and the
second is the threshold.

Only boxes containing \(m\) affect \(G_{\{m\}}\); in each box counted
by \(\sigma_n^\ast\), that index occurs once.  By symmetry, we may
put it first and multiply by \(k\).  To allow repeated indices
later, for every \(B=I_m\times I_{m_2}\times\cdots\times I_{m_k}\)
and \(0<u<1\), define the slice sum
\[
 F_B^{(m)}(u)
 =\sum_{x\in\{2m+u,\,2m+2-u\}}
  \int_{I_{m_2}\times\cdots\times I_{m_k}}
  \mathbf{1}_{\{x^k+\sum_{\nu=2}^k y_\nu^k\le n\}}
  \,dy_2\cdots dy_k.
\]
For distinct pair indices, independence gives
\(F_B^{(m)}(U_m)=\E[\delta_B\mid U_m]\).  For every box, including
those with repeated indices, \eqref{eq:pair-uniform} gives
\[
 \E F_B^{(m)}(U_m)
 =\int_B\mathbf{1}_{\{\sum_{\nu=1}^k y_\nu^k\le n\}}
       \,dy_1\cdots dy_k=d_B.
\]
Thus for distinct indices the law of total expectation also gives
\(d_B=\E[\E[\delta_B\mid U_m]]\).  Summing over \(m\), we obtain
\[
 \begin{aligned}
 \sum_m G_{\{m\}}
 &=\sum_m\bigl(\E[\sigma_n^\ast\mid U_m]-\E\sigma_n^\ast\bigr)\\
 &=\sum_m\sum_{\substack{B\text{ with distinct}\\
                         \text{pair indices}}}
       \bigl(\E[\delta_B\mid U_m]-d_B\bigr)\\
 &=k\sum_m\sum_{\substack{m_2,\ldots,m_k\ge0\\
                         m,m_2,\ldots,m_k\text{ distinct}}}
       \bigl(F_B^{(m)}(U_m)-\E F_B^{(m)}(U_m)\bigr).
 \end{aligned}
\]

Extend this centered slice sum to all \(m_2,\ldots,m_k\ge0\).
The added boxes have repeated indices.  Off the boundary,
\(F_B^{(m)}\) is constant, equal to \(2^k\) or \(0\), so its
centered contribution vanishes.  On a boundary box,
\(0\le F_B^{(m)}(U_m),\E F_B^{(m)}(U_m)\le2^k\), so the
centered contribution has absolute value at most \(2^k\).
Each added ordered box occurs once when we sum over its first
index \(m\).  Hence the total correction is at most
\[
 k2^k\#\{B:\ B\text{ boundary with a repeated index}\}
 \ll_k R^{k-2},
\]
by Lemma~\ref{lem:boundary}.

The remaining boxes now partition the full orthant, and their
slice integrals sum to
\[
 \begin{aligned}
 \sum_{m_2,\ldots,m_k\ge0}F_B^{(m)}(U_m)
 &=\sum_{s=2m}^{2m+1}\int_{[0,\infty)^{k-1}}
       \mathbf{1}_{\{\theta_s^k+\sum_{\nu=2}^k y_\nu^k\le n\}}
       \,dy_2\cdots dy_k\\
 &=f_n(\theta_{2m})+f_n(\theta_{2m+1})=q_{m,n}(U_m).
 \end{aligned}
\]
Subtracting expectations therefore yields
\begin{equation}\label{eq:linear-q}
 \sum_m G_{\{m\}}
 =k\sum_m\bigl(q_{m,n}(U_m)-\E q_{m,n}(U_m)\bigr)+O_k(R^{k-2}).
\end{equation}

It remains to bound the variation of \(q_{m,n}\).  Away from \(R\),
the two points move in opposite directions and the difference of
their first derivatives can be controlled by \(f_n''\).  Near \(R\)
we use a direct bound on the slice volumes.

First suppose \(2m+2\le R-1\), so both arguments of \(f_n\) stay
below \(R\).  We differentiate \(q_{m,n}\), then apply the mean value
theorem first to \(f_n'\) and then to \(q_{m,n}\).  For
\(u,u_1,u_2\in(0,1)\), this gives
\[
 \begin{aligned}
 q_{m,n}'(u)&=f_n'(2m+u)-f_n'(2m+2-u),\\
 |q_{m,n}'(u)|
 &\le\bigl|(2m+u)-(2m+2-u)\bigr|\sup_{x\in I_m}|f_n''(x)|\\
 &=2(1-u)\sup_{x\in I_m}|f_n''(x)|
 \le2\sup_{x\in I_m}|f_n''(x)|,\\
 |q_{m,n}(u_1)-q_{m,n}(u_2)|
 &\le2|u_1-u_2|\sup_{x\in I_m}|f_n''(x)|
 \le2\sup_{x\in I_m}|f_n''(x)|.
 \end{aligned}
\]
For \(0\le x<R\), put \(v=x/R\).  Using \(n=R^k\), the chain rule,
and \(1-v^k\ge1-v\), we obtain
\[
 \begin{aligned}
 f_n(x)&=C_{k-1,k}(R^k-x^k)^{1-1/k}
        =C_{k-1,k}R^{k-1}(1-v^k)^{1-1/k},\\
 f_n'(x)&=-C_{k-1,k}(k-1)R^{k-2}v^{k-1}(1-v^k)^{-1/k},\\
 |f_n''(x)|
 &=C_{k-1,k}(k-1)R^{k-3}
   \frac{v^{k-2}\bigl((k-1)-(k-2)v^k\bigr)}
        {(1-v^k)^{1+1/k}}\\
 &\ll_k R^{k-3}(1-v^k)^{-1-1/k}\\
 &\ll_k R^{k-3}(1-v)^{-1-1/k}
  =R^\gamma(R-x)^{-1-1/k}.
 \end{aligned}
\]
The mean \(\E q_{m,n}(U_m)\) lies between the infimum and supremum
of \(q_{m,n}\).  The distances \(R-(2m+2)\) in the following sum are
at least one and spaced by two, so the derivative bound and the
integral test give
\[
 \begin{aligned}
 \sum_{\substack{m\ge0\\2m+2\le R-1}}
       |q_{m,n}(U_m)-\E q_{m,n}(U_m)|
 &\le2\sum_{\substack{m\ge0\\2m+2\le R-1}}
       \sup_{x\in I_m}|f_n''(x)|\\
 &\ll_k R^\gamma
       \sum_{\substack{m\ge0\\2m+2\le R-1}}
       \bigl[R-(2m+2)\bigr]^{-1-1/k}\\
 &\ll_k R^\gamma\sum_{a=0}^{\infty}(1+2a)^{-1-1/k}\\
 &\ll_k R^\gamma\left(1+\int_0^\infty(1+2t)^{-1-1/k}\,dt\right)
   =\left(1+\frac{k}{2}\right)R^\gamma .
 \end{aligned}
\]
We now consider the intervals near \(R\).  Those starting at or above
\(R\) give zero, so the remaining indices satisfy
\[
 2m+2>R-1,\qquad 2m<R,
 \qquad\text{or equivalently}\qquad
 \frac{R-3}{2}<m<\frac R2.
\]
For these intervals, \(x\in I_m\) implies \(x\ge2m>R-3\).
Using the factorization of \(R^k-x^k\) and
\((k-1)(1-1/k)=\gamma\), we obtain
\[
 \begin{aligned}
 (n-x^k)_+
 &=(R^k-x^k)_+\le kR^{k-1}(R-x)_+\le3kR^{k-1},\\
 0\le f_n(x)
 &\le C_{k-1,k}(3kR^{k-1})^{1-1/k}\ll_k R^\gamma,\\
 0\le q_{m,n}(u)
 &=f_n(2m+u)+f_n(2m+2-u)\ll_k R^\gamma
 \qquad(0<u<1).
 \end{aligned}
\]
There are at most two such indices.  Since \(q_{m,n}\) is nonnegative,
\[
 \begin{aligned}
 \sum_{\substack{m\ge0\\(R-3)/2<m<R/2}}
       |q_{m,n}(U_m)-\E q_{m,n}(U_m)|
 &=\sum_{\substack{m\ge0\\(R-3)/2<m<R/2}}
       \left|q_{m,n}(U_m)-\int_0^1q_{m,n}(u)\,du\right|\\
 &\le\sum_{\substack{m\ge0\\(R-3)/2<m<R/2}}
       \sup_{0<u<1}q_{m,n}(u)\ll_k R^\gamma.
 \end{aligned}
\]

Combining both ranges with \eqref{eq:linear-q} and using
\(\gamma=k-2+1/k>k-2\), we obtain, for every choice of the points,
\begin{equation}\label{eq:first-order}
 \begin{aligned}
 \left|\sum_m G_{\{m\}}\right|
 &\le k\Biggl(\sum_{\substack{m\ge0\\2m+2\le R-1}}
       +\sum_{\substack{m\ge0\\(R-3)/2<m<R/2}}\Biggr)
       |q_{m,n}(U_m)-\E q_{m,n}(U_m)|+O_k(R^{k-2})\\
 &\ll_k R^\gamma+R^{k-2}\ll_k R^\gamma.
 \end{aligned}
\end{equation}

To show that \(\mathcal H_n=o(R^\gamma)\), we group the terms
\(G_T\) by their two largest indices and count boundary boxes
containing those indices.  The indices \(i,j\) used for these
groups are pair indices.  For \(T\subseteq S\) with
\(2\le|T|\le k\), put
\[
 M_n(T)=\sum_{\substack{
    j_1,\ldots,j_k\in S,\quad j_1<\cdots<j_k\\
    T\subseteq\{j_1,\ldots,j_k\}\\
    I_{j_1}\times\cdots\times I_{j_k}\text{ is a boundary box}
 }}1.
\]
Thus \(M_n(T)\) counts completions of \(T\) to a boundary box,
once per index set in increasing order.  For \(k=3\), a boundary
box with index set \(\{1,2,3\}\) contributes once to
\(M_n(\{1,2\})\).  This count is deterministic.

The square sum for the two prescribed indices is
\[
 D_2(n)=\sum_{\substack{T\subseteq S\\ |T|=2}}M_n(T)^2 .
\]
For fixed \(T\), the square counts ordered pairs of index sets
containing \(T\), including a set paired with itself.  They may
share further indices; coordinate orders are not counted.

We first relate \(D_2(n)\) to the moments of \(\mathcal H_n\).
Write \(U_{<i}=U_{S\cap[0,i)}\) and \(U_{\le i}=U_{S\cap[0,i]}\).
Every \(T\) with \(|T|\ge2\) is uniquely of the form
\(J\cup\{i,j\}\), where \(i<j\) are its two largest indices and
\(J\subseteq S\cap[0,i)\).  Hence
\[
 \begin{aligned}
 \mathcal H_n
 &=\sum_{j\in S}\sum_{\substack{\max T=j\\|T|\ge2}}G_T\\
 &=\sum_{j\in S}\sum_{\substack{i\in S\\i<j}}
       \sum_{J\subseteq S\cap[0,i)}G_{J\cup\{i,j\}} .
 \end{aligned}
\]
For example, if \(1,2,3\in S\), then \(G_{\{1,2,3\}}\) and
\(G_{\{2,3\}}\) occur in the same group \(i=2,j=3\).
The former corresponds to \(J=\{1\}\) and the latter to
\(J=\varnothing\).

To express each group through the boxes, apply
\eqref{eq:hoeffding-conditioning} with \(W=\sigma_n^\ast\).
For any fixed conditioning set \(Q\subseteq S\),
\[
 \E[\sigma_n^\ast\mid U_Q]=\sum_{T\subseteq Q}G_T .
\]
Here \(Q\) is fixed and \(T\) is summed over its subsets;
\(J\) denotes the remaining indices in a group below.
For \(i<j\), subtracting the sums with and without \(j\) gives
\[
 \E[\sigma_n^\ast\mid U_{\le i},U_j]
 -\E[\sigma_n^\ast\mid U_{\le i}]
 =\sum_{J\subseteq S\cap[0,i]}G_{J\cup\{j\}} .
\]
Subtract the same identity with \(U_{\le i}\) replaced by
\(U_{<i}\).  Only terms containing \(i\) remain, giving
\begin{equation}\label{eq:two-index-group}
 \begin{aligned}
 \sum_{J\subseteq S\cap[0,i)}G_{J\cup\{i,j\}}
 ={}&\E[\sigma_n^\ast\mid U_{\le i},U_j]
    -\E[\sigma_n^\ast\mid U_{\le i}]\\
    &-\E[\sigma_n^\ast\mid U_{<i},U_j]
    +\E[\sigma_n^\ast\mid U_{<i}].
 \end{aligned}
\end{equation}

We expand \eqref{eq:two-index-group} over the boxes counted by
\(\sigma_n^\ast\).  If a box omits \(r\in\{i,j\}\), independence gives
\(\E[\delta_B\mid U_Q,U_r]-\E[\delta_B\mid U_Q]=0\) for every
\(Q\subseteq S\setminus\{r\}\), so the terms that differ only by
conditioning on \(U_r\) cancel in pairs.  For a constant box count,
the four terms also sum to zero.  Thus in the sums below \(B\) runs
over the ordered boundary boxes with distinct indices containing
both \(i\) and \(j\).  Since \(0\le\delta_B\le2^k\), we have almost surely
\[
 \begin{aligned}
 \left|\sum_{J\subseteq S\cap[0,i)}G_{J\cup\{i,j\}}\right|
 &\le \sum_B\Bigl|
   \E[\delta_B\mid U_{\le i},U_j]-\E[\delta_B\mid U_{<i},U_j]
   -\bigl(\E[\delta_B\mid U_{\le i}]-\E[\delta_B\mid U_{<i}]\bigr)
 \Bigr|\\
 &\le 4\cdot2^k\sum_B1
   =4\cdot2^k k!M_n(\{i,j\}).
 \end{aligned}
\]
The last equality counts the \(k!\) coordinate orders of each
increasing tuple counted by \(M_n(\{i,j\})\).

For fixed \(j\), first reveal \(U_j\), then reveal the draws \(U_i\)
in increasing order for \(i<j\).  The group indexed by \((i,j)\)
depends only on \(U_{\le i},U_j\), the draws revealed so far,
and every term contains \(i\).  For the outer sums, reveal \(U_j\)
in increasing order.  The sum with largest index \(j\) depends
only on \(U_{\le j}\), and every term contains \(j\).
By independence and the conditional mean-zero property of \(G_T\),
\[
 \begin{aligned}
 \E\left[\sum_{J\subseteq S\cap[0,i)}G_{J\cup\{i,j\}}
       \,\middle|\,U_{<i},U_j\right]&=0 &&(i<j),\\
 \E\left[\sum_{\substack{T\subseteq S\\\max T=j,\ |T|\ge2}}G_T
       \,\middle|\,U_{<j}\right]&=0 &&(j\in S).
 \end{aligned}
\]
Thus each inner group has mean zero given the earlier draws and
\(U_j\), while each outer sum has mean zero given the earlier
draws.  Both sequences are therefore martingale differences.

For \(p\ge2\), apply Lemma~\ref{lem:martingale} first to the sums
indexed by \(j\), and then to the groups indexed by \(i\) for each
fixed \(j\).  The box bound gives the last inequality below.
\[
 \begin{aligned}
 \|\mathcal H_n\|_p^2
 &\le (p-1)\sum_{j\in S}
      \left\|\sum_{\substack{\max T=j\\|T|\ge2}}G_T\right\|_p^2\\
 &=(p-1)\sum_{j\in S}
      \left\|\sum_{\substack{i\in S\\i<j}}
       \sum_{J\subseteq S\cap[0,i)}G_{J\cup\{i,j\}}\right\|_p^2\\
 &\le (p-1)^2\sum_{j\in S}\sum_{\substack{i\in S\\i<j}}
      \left\|\sum_{J\subseteq S\cap[0,i)}G_{J\cup\{i,j\}}\right\|_p^2\\
 &\le (4\cdot2^k k!)^2(p-1)^2D_2(n).
 \end{aligned}
\]
Taking square roots and setting \(p=2\ell\), we obtain
\begin{equation}\label{eq:moment-H}
 \|\mathcal H_n\|_{2\ell}\ll_k(2\ell-1)D_2(n)^{1/2}.
\end{equation}
Choose \(c>0\), depending only on \(k\), such that the bound in
\eqref{eq:moment-H} holds with constant \(c\).  For \(n\ge2\), set
\(\ell=\lceil2\log n\rceil\), so \(2\ell-1\le4\log n+1\).
Markov's inequality gives, when \(D_2(n)>0\),
\[
 \begin{aligned}
 &\Prob\bigl(|\mathcal H_n|>e\cdot c(4\log n+1)D_2(n)^{1/2}\bigr)\\
 &\qquad\le
 \frac{\E|\mathcal H_n|^{2\ell}}
      {\bigl(e\cdot c(4\log n+1)D_2(n)^{1/2}\bigr)^{2\ell}}
 \le e^{-2\ell}\le n^{-4}.
 \end{aligned}
\]
If \(D_2(n)=0\), then \(\mathcal H_n=0\) almost surely, so the
same tail bound holds.  Define the exceptional events
\[
 \mathcal E_n=\left\{
   |\mathcal H_n|>e\cdot c(4\log n+1)D_2(n)^{1/2}
 \right\},\qquad n\ge2.
\]
Each \(\mathcal E_n\) is measurable on the common probability
space of \(U_0,U_1,\ldots\): \(\mathcal H_n\) is a measurable
function of finitely many draws, and \(D_2(n)\) is deterministic.
The preceding tail estimate verifies the summability condition:
\[
 \sum_{n=2}^{\infty}\Prob(\mathcal E_n)
 \le\sum_{n=2}^{\infty}n^{-4}<\infty.
\]
Thus Lemma~\ref{lem:borel-cantelli} applies, and with probability
one only finitely many \(\mathcal E_n\) occur.  Consequently,
\begin{equation}\label{eq:simultaneous}
 |\mathcal H_n|
 \le e\cdot c(4\log n+1)D_2(n)^{1/2}
 \ll_k(\log n)D_2(n)^{1/2}
\end{equation}
for all sufficiently large integers \(n\), using
\(4\log n+1\ll\log n\) for \(n\ge2\).
The starting index may depend on the realization, while the
implicit constant depends only on \(k\).  No independence between
the events \(\mathcal E_n\) is needed.

By \eqref{eq:simultaneous}, a power saving over \(R^{2\gamma}\)
in \(D_2(n)\) suffices.  We first bound \(M_n(T)\) for fixed \(T\).
For a boundary box, each \(2m_\nu\le R\).  Since
\((x+2)^k-x^k\) increases with \(x\ge0\), putting
\(L=k((R+2)^k-R^k)\asymp_k R^{k-1}\) gives
\[
 \begin{aligned}
 0&\le\sum_{\nu=1}^k\bigl((2m_\nu+2)^k-(2m_\nu)^k\bigr)\le L,\\
 (n-L)_+&\le\sum_{\nu=1}^k y_\nu^k\le n+L
 \qquad(y\in B).
 \end{aligned}
\]
For fixed \(T=\{m_1,m_2\}\), \(m_1<m_2\), with \(M_n(T)>0\), put
\[
 t=n-(2m_1)^k-(2m_2)^k\ge0.
\]
In the sums below, \(t\) always takes the value corresponding to
\(T\).  Place the two prescribed coordinates first and fix them
at \(2m_1,2m_2\).  Relabeling the remaining coordinates as
\(y_1,\ldots,y_{k-2}\), their boxes lie in the layer
\[
 (t-L)_+\le\sum_{\nu=1}^{k-2}y_\nu^k\le t+L.
\]
In increasing order of the remaining indices, the boxes have
disjoint interiors and volume \(2^{k-2}\).  The volume formula
\eqref{eq:volume} therefore gives
\[
 \begin{aligned}
 M_n(T)
 &\le\frac{C_{k-2,k}}{2^{k-2}}
       \bigl((t+L)^{1-2/k}-(t-L)_+^{1-2/k}\bigr)\\
 &\ll_k
 \begin{cases}
  L(t-L)^{-2/k}\ll_k Lt^{-2/k},&t\ge2L,\\
  (t+L)^{1-2/k}\ll_k L^{1-2/k},&0\le t<2L.
 \end{cases}
 \end{aligned}
\]
The first case uses the mean value theorem and \(t-L\ge t/2\);
in the second, \(L\le t+L<3L\).  Since \(t+L\le3t/2\) in the
first case, both estimates give
\begin{equation}\label{eq:completion-bound}
 M_n(T)\ll_k\frac{L}{(t+L)^{2/k}}.
\end{equation}

To bound the total number of boundary boxes, use the full
\(k\)-dimensional layer above.  The boxes have disjoint interiors
and volume \(2^k\), so \eqref{eq:volume} with \(d=k\) gives
\[
 2^k\sum_{B\text{ boundary}}1
 \le C_k\bigl((n+L)-(n-L)_+\bigr)
 \le2C_kL\ll_k R^{k-1}.
\]
Each boundary box \(B=I_{j_1}\times\cdots\times I_{j_k}\) with
\(j_1<\cdots<j_k\) contains \(\binom{k}{2}\) two-element sets,
so double counting gives
\[
 \sum_{\substack{T\subseteq S\\|T|=2}}M_n(T)
 =\binom{k}{2}\sum_{\substack{B\text{ boundary}\\j_1<\cdots<j_k}}1
 \ll_k R^{k-1}.
\]

To estimate \(D_2(n)\), omit zero terms and split at
\(t=LR^{1/4}\).  In the first range, the rectangles
\(I_{m_1}\times I_{m_2}\) lie in the layer
\[
 n-LR^{1/4}\le (2m_1)^k+(2m_2)^k
 \le x_1^k+x_2^k\le n+L .
\]
Since \(LR^{1/4}/n\ll_k R^{-3/4}\), the two thresholds are comparable
to \(n\) for large \(R\).  The rectangles have disjoint interiors
and area four.  We bound their union by the full layer; the
restriction \(m_1<m_2\) only reduces the area being counted.
Thus \eqref{eq:volume} and the mean value theorem give
\[
 \begin{aligned}
 \sum_{\substack{T\subseteq S,\ |T|=2\\M_n(T)>0,\ t\le LR^{1/4}}}1
 &\le\frac{C_{2,k}}4
       \bigl((n+L)^{2/k}-(n-LR^{1/4})^{2/k}\bigr)\\
 &\ll_k LR^{1/4}n^{2/k-1}\ll_k R^{5/4}.
 \end{aligned}
\]
In the first range, \eqref{eq:completion-bound} gives the uniform
bound \(L^{1-2/k}\).  In the second, \(t+L>LR^{1/4}\), so
\[
 M_n(T)\ll_k\frac{L}{(LR^{1/4})^{2/k}}
 =L^{1-2/k}R^{-1/(2k)}.
\]
Hence
\[
 M_n(T)\ll_k
 \begin{cases}
  L^{1-2/k},&0\le t\le LR^{1/4},\\
  L^{1-2/k}R^{-1/(2k)},&t>LR^{1/4}.
 \end{cases}
\]
In the first range, use the number of rectangles and the squared
bound on \(M_n(T)\).  In the second, bound one factor and sum the
other by the double count above.  Since
\((3/4-2/k)-1/(4k)=3(k-3)/(4k)\ge0\) for \(k\ge3\), this gives
\begin{equation}\label{eq:codegrees}
 \begin{aligned}
 D_2(n)
 &=\sum_{\substack{|T|=2,\ M_n(T)>0\\t\le LR^{1/4}}}M_n(T)^2
   +\sum_{\substack{|T|=2,\ M_n(T)>0\\t>LR^{1/4}}}M_n(T)^2\\
 &\ll_k R^{5/4}L^{2-4/k}
       +L^{1-2/k}R^{-1/(2k)}\sum_{|T|=2}M_n(T)\\
 &\ll_k R^{2\gamma-(3/4-2/k)}+R^{2\gamma-1/(2k)}
 \ll_k R^{2\gamma-1/(4k)}.
 \end{aligned}
\end{equation}

By \eqref{eq:codegrees} and \eqref{eq:simultaneous}, almost surely,
\[
 \frac{|\mathcal H_n|}{R^\gamma}
 \ll_k R^{-1/(8k)}\log n
 =n^{-1/(8k^2)}\log n\longrightarrow0,
\]
which proves \(\mathcal H_n=o(R^\gamma)\).

Returning to \eqref{eq:total-discrepancy}, we combine
\eqref{eq:first-order}, \(\mathcal H_n=o(R^\gamma)\), and
\(\gamma>k-2\) to obtain, almost surely,
\begin{equation}\label{eq:continuous}
 \sigma_n=C_kn+O_k(R^\gamma)
     =C_kn+O_k\bigl(n^{1-2/k+1/k^2}\bigr)
\end{equation}
for every sufficiently large integer \(n\).

Choose a realization for which \eqref{eq:continuous} holds for all
sufficiently large \(n\).  Shifting \(n\) by \(k\) changes the main
term by only \(kC_k\), so \(\sigma_{n-k}=C_kn+O_k(R^\gamma)\) as well.
The rounding comparison \eqref{eq:rounding} therefore gives, for
some \(c_k>0\) and all sufficiently large \(n\),
\[
 C_kn-c_kR^\gamma
 \le \sigma_{n-k}
 \le \sum_{h\le n}r_{kA}(h)
 \le \sigma_n
 \le C_kn+c_kR^\gamma.
\]
This proves \eqref{eq:main}, since
\(C_k=\Gamma(1+1/k)^k\) and \(R^\gamma=n^{1-2/k+1/k^2}\).

Finally, the bounds \(s^k+1\le a_s\le(s+1)^k\) give, for \(X\ge1\),
\[
 s+1\le X^{1/k}\ \Longrightarrow\ a_s\le X,
 \qquad s\ge X^{1/k}\ \Longrightarrow\ a_s>X.
\]
Hence \(\lfloor X^{1/k}\rfloor\le A(X)\le\lceil X^{1/k}\rceil\),
giving the asserted density.
\qed

\section{Fluctuations in the paired model}\label{sec:fluctuations}

The proof of Theorem~\ref{thm:main} shows that individual pairs give
the largest error bound.  We now show that a pair near \(R\) can
produce fluctuations of that order for infinitely many \(n\).

\begin{proposition}\label{prop:fluctuation}
Fix \(k\ge3\).  Let \(A\) be obtained from the points in
\eqref{eq:points} by taking powers and rounding as above.
With probability one,
\[
 \limsup_{n\to\infty}
 \frac{|\sum_{h\le n}r_{kA}(h)-C_kn|}{n^{1-2/k+1/k^2}}>0.
\]
\end{proposition}

\begin{proof}
Retain \(R=n^{1/k}\) and \(\gamma=k-2+1/k\).
By \eqref{eq:total-discrepancy}, \eqref{eq:linear-q}, and
\(\mathcal H_n=o(R^\gamma)\), almost surely,
\[
 \sigma_n-C_kn=k\sum_{i\in S}\bigl(q_{i,n}(U_i)-\E q_{i,n}(U_i)\bigr)+o(R^\gamma)
\]
as \(n\) tends to infinity through the integers.

For each \(m\in\Nzero\), take \(R=2m+1\) and \(n=R^k\).  Then
\(S=\{0,\ldots,m\}\), so the centered sum above depends only on
\(U_0,\ldots,U_m\).  In pair \(m\), only the
lower point lies below \(R\).  By the mean value theorem,
\(R^k-(R-1+u)^k=kR^{k-1}(1-u)(1+O_k(R^{-1}))\) uniformly for
\(0<u<1\).  Since \((k-1)(1-1/k)=\gamma\), we have
\[
 \begin{aligned}
 q_{m,n}(u)&=C_{k-1,k}\bigl(R^k-(R-1+u)^k\bigr)^{1-1/k}\\
 &=C_{k-1,k}k^{1-1/k}R^\gamma(1-u)^{1-1/k}
     \bigl(1+O_k(R^{-1})\bigr).
 \end{aligned}
\]
Evaluating at \(u=1/4,3/4\) gives
\[
 \begin{aligned}
 k\bigl(q_{m,n}(1/4)-q_{m,n}(3/4)\bigr)
 &=C_{k-1,k}k^{2-1/k}R^\gamma
   \bigl((3/4)^{1-1/k}-(1/4)^{1-1/k}+O_k(R^{-1})\bigr)\\
 &\ge2b_kR^\gamma
 \end{aligned}
\]
for some \(b_k>0\) and all sufficiently large \(m\).

Fix \(U_0,\ldots,U_{m-1}\).  The centered sum, including the
factor \(k\), is a decreasing function of \(U_m\).  Its values
at \(1/4\) and \(3/4\) differ by at least
\(2b_kR^\gamma\).

If the value at \(1/4\) is at least \(b_kR^\gamma\), so is every
value on \((0,1/4)\).  Otherwise every value on \((3/4,1)\) is at
most \(-b_kR^\gamma\).  Thus, given \(U_0,\ldots,U_{m-1}\), the
absolute value of the centered sum is at least \(b_kR^\gamma\)
with probability at least \(1/4\).

For \(R=2m+1\) and \(n=(2m+1)^k\), define \(E_m\) as below.
Thus every \(q_{i,n}\) in \(E_m\) uses this same threshold, which
changes with \(m\).  The preceding argument and successive
conditioning give, for any sufficiently large \(m_0\) and integer
\(N\ge1\),
\[
 \begin{gathered}
 E_m=\left\{\left|k\sum_{i=0}^{m}
       \bigl(q_{i,n}(U_i)-\E q_{i,n}(U_i)\bigr)\right|<b_kR^\gamma\right\},\\
 \Prob(E_m\mid U_0,\ldots,U_{m-1})\le\frac34,
 \end{gathered}
\]
\[
 \begin{aligned}
 \Prob\left(\bigcap_{m=m_0}^{m_0+N-1}E_m\right)
 &=\E\left[\mathbf{1}_{\bigcap_{m=m_0}^{m_0+N-2}E_m}\,
       \Prob(E_{m_0+N-1}\mid U_0,\ldots,U_{m_0+N-2})\right]\\
 &\le\frac34\Prob\left(\bigcap_{m=m_0}^{m_0+N-2}E_m\right)
 \le\left(\frac34\right)^N,\\
 \Prob\left(\bigcap_{m\ge m_0}E_m\right)
 &=\lim_{N\to\infty}\Prob\left(\bigcap_{m=m_0}^{m_0+N-1}E_m\right)=0.
 \end{aligned}
\]
Taking the countable union over \(m_0\) shows that \(E_m\) fails
infinitely often almost surely.  The \(o(R^\gamma)\) remainder then
gives, for infinitely many \(n=(2m+1)^k\),
\[
 |\sigma_n-C_kn|\ge b_kR^\gamma-o(R^\gamma)
                   \ge\frac{b_k}{2}R^\gamma .
\]

For these \(n\), apply \eqref{eq:rounding} at \(n\) in the negative
case and at \(n+k\) in the positive case.
\[
 \begin{aligned}
 \sum_{h\le n}r_{kA}(h)-C_kn
 &\le \sigma_n-C_kn\le-\tfrac{b_k}{2}R^\gamma
 &&\text{if \(\sigma_n-C_kn<0\)},\\
 \sum_{h\le n+k}r_{kA}(h)-C_k(n+k)
 &\ge \sigma_n-C_kn-kC_k\ge\tfrac{b_k}{2}R^\gamma-kC_k
 &&\text{if \(\sigma_n-C_kn>0\)}.
 \end{aligned}
\]
Since \(R^\gamma=n^{1-2/k+1/k^2}\) and \((n+k)/n\to1\), these bounds
imply, with probability one,
\[
 \limsup_{n\to\infty}
 \frac{\left|\sum_{h\le n}r_{kA}(h)-C_kn\right|}
      {n^{1-2/k+1/k^2}}
 \ge\frac{b_k}{2}>0.
\]
\end{proof}

\end{document}